\documentclass
{amsart}
\usepackage{amssymb,mathrsfs}
\usepackage{amssymb}
\usepackage{amsmath}
\usepackage{amsfonts, mathrsfs}
\usepackage{mathtools}
\usepackage{epsfig}
\usepackage{color}
\usepackage{epstopdf}
\usepackage{graphicx}
\usepackage{srcltx}
\usepackage{amsthm}
\usepackage[mathscr]{eucal}
\usepackage{enumerate}
\usepackage{verbatim}
\usepackage[draft]{todonotes}
\usepackage{relsize}
\usepackage{hyperref}
\usepackage{cases} 
\usepackage{empheq}
\usepackage{subfigure}
\usepackage{tikz}
\usepackage{array}
\usepackage{tabularx}
\newcolumntype{C}[1]{>{\centering\arraybackslash}p{#1}}

\usetikzlibrary{arrows,automata}
\usepackage{varwidth}
\usepackage{tikz}
\usetikzlibrary{backgrounds,patterns,calc}
\usepackage{xparse}

\usepackage{enumitem}

\usepackage[normalem]{ulem}
\newcommand{\stkout}[1]{\ifmmode\text{\sout{\ensuremath{#1}}}\else\sout{#1}\fi}

\hypersetup{
	colorlinks=false,       
	linkcolor=blue,          
	citecolor=red,        
	filecolor=magenta,      
	urlcolor=cyan           
}

\newcommand{\Be}{\begin{equation}}
\newcommand{\Ee}{\end{equation}}
\newcommand{\Bea}{\begin{eqnarray}}
\newcommand{\Eea}{\end{eqnarray}}
\newcommand{\Bel}{\begin{align}}
\newcommand{\Eel}{\end{align}}
\newcommand{\Beas}{\begin{eqnarray*}}
	\newcommand{\Eeas}{\end{eqnarray*}}
\newcommand{\Benu}{\begin{enumerate}}
	\newcommand{\Eenu}{\end{enumerate}}
\newcommand{\Bi}{\begin{itemize}}
	\newcommand{\Ei}{\end{itemize}}

\newcommand{\B}{\Big}
\renewcommand{\S}{\mathbb S}
\numberwithin{equation}{section}
\newcommand{\supp} {\text{supp\! }}

\theoremstyle{plain}
\newtheorem{thm}{Theorem}[section]

\newtheorem{lem}[thm]{Lemma}

\newtheorem{conj}[thm]{Conjecture}

\theoremstyle{remark}

\theoremstyle{definition}

\definecolor{ao}{rgb}{0, 0.5, 0}

\newcommand{\wa}{e^{it\sqrt{\mathsmaller{-}\Delta}}}
 \newcommand{\pl}{P_\lambda}
 
 \newcommand{\Bb}{\mathbb B}

\usetikzlibrary{patterns}
\tikzset{
    slope/.store in=\slope,
    pattern color/.append code={\def\patterncolor{#1}}
}

\newcommand{\theslope}{0.7}

\pgfdeclarepatternformonly{diagonal lines}
{\pgfpoint{-.1mm/\theslope}{-.1mm}}{\pgfpoint{1.1mm/\theslope}{1.1mm}}
{\pgfpoint{1mm/\theslope}{1mm}}
{
    \pgfsetlinewidth{0.4pt}
    \pgfpathmoveto{\pgfpoint{-.1mm/\theslope}{-.1mm}}
    \pgfpathlineto{\pgfpoint{1.1mm/\theslope}{1.1mm}}
    \pgfusepath{stroke}
}

\pgfdeclarepatternformonly[\slope,\patterncolor]{slant lines}
{\pgfpoint{-.1mm/\slope}{-.1mm}}{\pgfpoint{1.1mm/\slope}{1.1mm}}
{\pgfpoint{1mm/\slope}{1mm}}
{
    \pgfsetlinewidth{0.4pt}
    \pgfpathmoveto{\pgfpoint{-.1mm/\slope}{-.1mm}}
    \pgfpathlineto{\pgfpoint{1.1mm/\slope}{1.1mm}}
    \pgfsetstrokecolor{\patterncolor}
    \pgfusepath{stroke}
}

\usepackage{tikz}
\usepackage{pgfplots}
\pgfplotsset{compat=1.13,
tick label style={color=white},
  label style={font=\small},
  legend style={font=\small}}
\usepackage{mathtools}

\title[divergence set of the wave equation]
{Dimension of divergence set  of \\ the wave equation}

\author{Seheon Ham}
\author{Hyerim Ko}
\author{Sanghyuk Lee}

\thanks{}
\address{Department of Mathematical Sciences and RIM, Seoul National University, Seoul 08826, Republic of Korea}
\email{seheonham@snu.ac.kr}
\email{kohr@snu.ac.kr}
\email{shklee@snu.ac.kr}

\begin{document}

\begin{abstract} 
We consider the Hausdorff dimension of the divergence set on which the pointwise convergence $\lim_{t\rightarrow 0} \wa f(x) = f(x)$ fails when  $f \in H^s(\mathbb R^d)$.
We especially prove the conjecture raised by Barcel\'o, Bennett, Carbery and Rogers \cite{BBCR} for $d=3$, and improve the previous results in higher dimensions $d\ge4$. 
We also show that  a Strichartz type estimate for $f\to \wa f$ with the measure $ dt\,d\mu(x)$ is essentially equivalent to  the estimate for the spherical average of  $\widehat \mu$ which has been extensively studied for the Falconer distance set problem. The equivalence provides shortcuts
to the recent results due to Liu \cite{Liu} and Rogers \cite{Rogers}.  
\end{abstract}

\maketitle

\section{introduction}\label{Sec:1}

Let $d\ge 1$. We consider  the wave operator 
\[
\wa f(x) = \frac{1}{(2\pi)^d} \int_{\mathbb R^d} e^{i (x\cdot \xi + t |\xi|)} \widehat f(\xi) d\xi.
\]
In this note we  are mainly concerned with the pointwise behavior of  $\wa f$ as $t\to 0$,  from  which we can deduce the 
pointwise convergence to the initial data  of the solution $u$ to  the Cauchy problem:
\[
\partial_t^2 u -  \Delta u =0  , \quad  u(\cdot,0)= u_0, \quad \partial_t u(\cdot,0)=  u_1,  
\]
where $(u_0,u_1) \in  H^s(\mathbb R^d)  \times  H^{s-1}(\mathbb R^d) $.
Here, $ H^s(\mathbb R^d)$ denotes the inhomogeneous Sobolev space of order $s$ which is  equipped with the norm $\|f\|_{ H^s} = \| (1+|\cdot|^2)^{\frac s 2} \widehat f\|_{L^2} $.

It is well-known that $\wa f$ converges to $f$  almost everywhere as $t\rightarrow 0$ for all $f \in H^s(\mathbb R^d)$ if $s>1/2$ while the convergence fails in general if $s\le \frac12$.
The convergence follows from the maximal estimate 
\begin{equation}\label{maximal} 
\| \sup_{0<t<1} |\wa f |\|_{L^2(\mathbb R^d)}  \lesssim \|f\|_{ H^s(\mathbb R^d) }
\end{equation} 
for $s>1/2$, which is  an  easy consequence of the Sobolev imbedding and Plancherel's theorem. (See for exmaple \cite{Cowling}.) 
The estimate fails for $s\le 1/2$ (see \cite{Walther}) and almost everywhere convergenece also fails for $s\le 1/2$. This can be shown using Stein's maximal theorem \cite{Stein2},  by which almost everywhere convergence for all $f\in H^{s}$ implies $H^{s}$--$L^{2,\infty}$ bound on the maximal operator $\sup_{0<t<1} |\wa f |$. 
 On the other hand,  it is not difficult to see  the $H^{1/2}$--$L^{2,\infty}$ estimate fails by a construction (see Lemma \ref{weak-type} in Appendix).

We study  a more refined pointwise behavior
of $\wa f$ by considering the Hausdorff dimension of the divergence set 
\[  
\mathfrak D_d  (f)  : = \{x:   \wa f(x)    \not\to f(x) ~\text{as}~ t\rightarrow 0 \}. 
\]
Dimension of $\mathfrak D_d  (f)$  naturally relates to  the regularity exponent  $s$ for which $f\in H^s$.  Almost everywhere convergence only tells  the divergence set is of Lebesgue measure zero. 
However, if the initial datum  $f$ possesses  an additional regularity, one may expect that the set $\mathfrak D_d$ has 
 smaller dimension. A natural question is how big the (Hausdorff) dimension of  $\mathfrak D_d$ can be depending the regularity of $f$. 

To study the question,  let  us set 
\begin{equation*}  
\mathcal D_d( s )  := \sup_{ f \in H^s(\mathbb R^d)}    \dim_H    \mathfrak D_d (f) . 
\end{equation*} 
The problem of determining $\mathcal D_d(s)$ was first considered by Barcel\'o, Bennett, Carbery and Rogers \cite{BBCR} 
for the more general dispersive equations and the result was later extended  by Luc\`{a} and Rogers \cite{LR}. 
Since $|\wa f |  \lesssim \|f\|_{H^s}$ for $s>d/2$, it follows 
 $\mathcal D_d(s) =0$ for $s >d/2$.  On the other hand, $\mathcal D_d(s)=d$ if $s\le 1/2$ since  the  $H^{1/2}$--$L^{2,\infty}$  bound   \eqref{maximal} generally fails.
 Thus it is sufficient to consider 
 \[1/2<  s\le d/2.\]
So,  there is nothing interesting if $d=1$.  For $d\ge 2$, the problem is no longer trivial.  When $d=2$ Barcel\'o et. al. \cite{BBCR}  obtained  the optimal result:    $\mathcal D_2(s) =2-2s$ for $3/4<s\le 1$ and  $\mathcal D_2(s) =4-4s$ for $1/2<s\le 3/4$. In higher dimensions the following was  conjectured to be true ({\cite[Theorem 5.1]{BBCR}}). 
 
\begin{conj} 
\label{conjecture}
For $d\ge 3$,  
\begin{align}\label{Dds}
\mathcal D_d(s) & = 
\begin{cases}
\, d-2s, \quad & 1\le s \le \frac d2,\\
 \, d+2-4s, \quad & \frac 12<s \le 1.
\end{cases}  
\end{align}
\end{conj}
\vspace{-6pt}
The lower bound was verified by considering specific functions and measures which show the sharpness of maximal estimate against the measures.  As is typical with the problems of similar nature, proving the upper bound is more difficult.

Let $\mathbb B^d(x,\rho)\subset \mathbb R^d$ denote the closed ball centered at $x$ of radius $\rho$.  
For $\alpha \in(0,d]$, we denote by $\mathfrak  M^d(\alpha)$ the collection of
non-negative Borel measures $\mu$ supported in $\Bb^d(0,1)$ such that 
\[
\mu ( \mathbb B^d (x, \rho)) \le C_\mu \rho^\alpha
\]
for all $x\in \mathbb R^d$ and $\rho>0$ with a constant $C_\mu$. 
If $\mu \in \mathfrak M^d(\alpha)$, we define 
\[
\langle \mu\rangle_\alpha = \sup_{x\in \mathbb R^{d},\, \rho>0} \rho^{-\alpha} 
\mu(  \mathbb B^d (x,\rho) ).
\]
An  approach to obtain the upper bound on $\mathcal D_d(s)$ is to exploit  the decay order $\beta=\beta(\alpha)$ of the estimate for $L^2$ average of $\widehat \mu$ over the sphere\,:
\begin{equation}\label{scale-beta2}     
\int_{\S^{d-1}} | \widehat {\mu}(\lambda\theta) |^2 d\theta \lesssim   \lambda^{-\beta} \|\mu\| \langle\mu\rangle_\alpha, \quad \mu\in  \mathfrak M^d(\alpha),
\end{equation}
where $\| \mu\| $ is the total variation of $\mu$.

The estimate  \eqref{scale-beta2} was previously  studied by various authors (\cite{Wolff2, Erdogan,  Erdogan2}) being related to the Falconer distance set problem, and further improvements were recently  obtained  
(\cite{LR, DZ}).  These improvements rely on sophisticated argument which has been  developed in the study of Fourier restriction problem.

An upper bound on $\mathcal D_d(s)$ can be obtained using \eqref{scale-beta2}. In fact, it was shown in \cite{BBCR} that  $\mathcal D_d(s) \le \alpha$ if  \eqref{scale-beta2} holds with $\beta >d-2s$. 
Thus, making use of the best estimate  so far (see \cite{Mattila, DZ}), 
 we have,  for $d\ge  3$, 
\begin{align}\label{div-set}
\mathcal D_d(s) \le 
\begin{cases}
\, d-2s, ~ & ~   \hfill \frac{d+1}{4} <  s \le \frac{d}2, \\
\, \frac{d^2-2ds } {d-1}, ~ & ~   \hfill \frac 12<s \le \frac{d+1}4. 
\end{cases} 
\end{align}
The bound in the case of $\frac{d+1}{4} <  s \le \frac{d}2$ coincides  with the sharp upper bound  of $\beta= \alpha$ for $0<\alpha\le \frac{d-1}{2}$ (see \cite{Mattila}).  
However, when $\frac12 <s \le \frac{d+1}{4}$ and $d\ge 4$, the current approach based on  \eqref{scale-beta2} does not seem to be efficient enough to  prove the optimal upper bound for $\mathcal D_d(s)$. 
In fact,  considering  the known  upper bound on $\beta$ (for example see \cite{Du}) for which \eqref{scale-beta2} holds, 
 the conjectured optimal upper bound $\mathcal D_d(s)=d+2-4s$ which corresponds to the estimate   \eqref{scale-beta2} with  $\beta = \frac{d+\alpha-2}{2}$ can not be achieved by the aforementioned implication from \eqref{scale-beta2}.

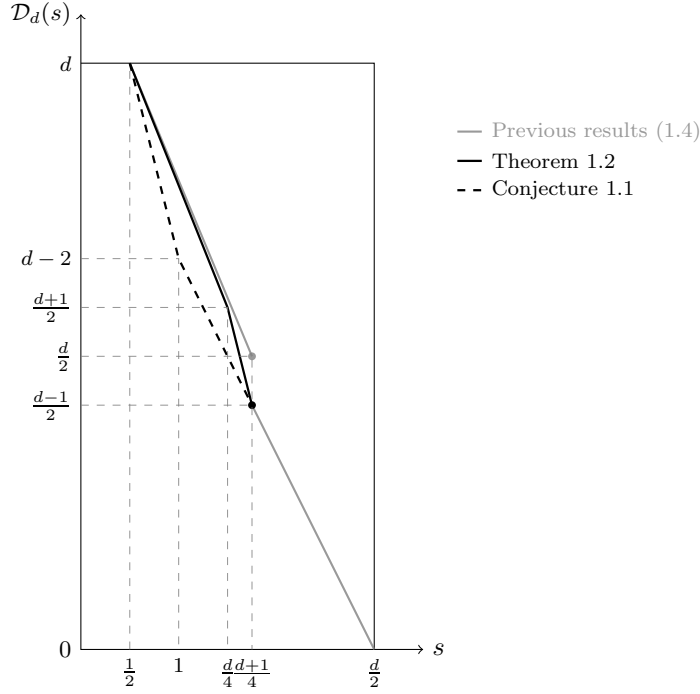
\begin{figure} 
	\begin{tikzpicture}[scale=1.3]
	\draw [<->] (0,6.5) node (yaxis) [left] {\small $\mathcal D_d(s)$}|- (3.5,0) node ( xaxis) [right] {$s$}  ;
	\draw (0,6)--(3,6)--(3,0);

	\draw [gray, thick, opacity=.8] (1.75, 2.5)--(3,0);
	\draw [gray, opacity=.8, thick] (0.5,6) --(1.75,3);
	\draw [line width=.3mm] (0.5,6)--(1.5, 3.5)--(1.75,2.5);
	\draw [dashed, line width=.3mm] (0.5,6)--(1,4)--(1.75,2.5);

 	\fill (1.75,2.5) circle (.04cm);
	\fill [gray, opacity=.8] (1.75,3) circle (.04cm);
	
	\draw (0.5,0) node [below] {\footnotesize $\frac 12$};
	\draw (1.5,0) node [below] {\footnotesize $\frac d4$};
	\draw (1.75,0) node [below] {\footnotesize $\frac {d+1}4$};
	\draw (1,0) node [below] {\footnotesize $1$};
	\draw (3,0) node [below] {\footnotesize $\frac d2$};
	
 	\draw [gray, dashed, opacity=.8] (0.5,0)--(0.5,6);
	\draw [gray, dashed, opacity=.8] (0,4)--(1,4)--(1,0);
	\draw [gray, dashed, opacity=.8] (0,3.5)--(1.5,3.5)--(1.5,0);
	\draw [gray, dashed, opacity=.8] (0,3)--(1.75,3)--(1.75,0);
	\draw [gray, dashed, opacity=.8] (0,2.5)--(1.75,2.5);

	\draw (0,0) node [left] {\small $0$};
	\draw (0,2.5) node [left] {\small $\frac {d-1}2$};
	\draw (0,3) node [left] {\small $\frac d2$};
	\draw (0,3.5) node [left] {\small $\frac {d+1}2$};
	\draw (0,4) node [left] {\footnotesize $d-2$};
	\draw (0,6) node [left] {\footnotesize $d$};
	
	\draw [gray, thick, opacity=.8] (3.85,5.3)--(4.1,5.3) node [right] {\footnotesize Previous results \eqref{div-set}};
	\draw [line width=.3mm] (3.85,5)--(4.1,5) node [right] {\footnotesize Theorem \ref{divergence}};
	\draw [dashed, line width=.3mm] (3.85,4.7)--(4.1,4.7) node [right] {\footnotesize Conjecture \ref{conjecture}};
	\end{tikzpicture}
	\caption{Dimension of the divergence set, $d\ge 6$} \label{div-pic}
\end{figure}

In this paper we take an alternative approach which relies on the fractal Strichartz estimate with respect to a measure (see \eqref{fracS} below), which  was previously studied by some authors (see \cite{Wolff3, Erdogan, CHL, Harris2, Rogers}).   
Via the approach we prove  the conjecture \eqref{Dds} when $d=3$ and  improve the previously known results (see \eqref{div-set}) for higher dimensions  $d\ge4$ and $\frac 12<s \le \frac{d+1}4$. 
\begin{thm}\label{divergence}
Let  $d \ge 3$. Then 
\begin{align*}
\mathcal D_d(s) \le
\begin{cases}
\ d-2s,  & \quad \frac {d+1}4 \le   s  \le   \frac{d}2,  
\\
\ \frac{3d+1}{2}-4s ,      & \quad \frac {d}4 \le  s \le  \frac{d+1}4, 
\\
\ \frac{d^2-d-1}{d-2}-\frac{2(d-1)s}{d-2},     &\quad 
 \hfill \frac12 < s \le 
\frac{d}4.  
\end{cases} 
\end{align*}
\end{thm}
In Figure \ref{div-pic} we compare the result in Theorem \ref{divergence},  the  conjectured optimal bound \eqref{Dds} and the previously known  result  \eqref{div-set}.

\subsection{Maximal estimates with respect to general measures.} 
To examine the size of the divergence set,
we consider the maximal estimate  
\begin{equation}\label{suptR}
\| \sup_{0<t <1}  |\wa  f   |  \|_{L^2 (d\mu)} \lesssim  \langle \mu \rangle_\alpha^{\frac12} \|f\|_{H^s(\mathbb R^d)}.
\end{equation}
It was shown in \cite{BBCR} that 
the  estimate \eqref{suptR} holds only if 
$ s\ge \max ( \frac{d-\alpha}{2},~  \frac{d+2-\alpha}{4})$. 
(See also Lemma \ref{nece} in Appendix A).  
In fact, the lower bound of $\mathcal D_d(s)$ was obtained by proving the failure of   \eqref{suptR} for $s< \max ( \frac{d-\alpha}{2},~  \frac{d+2-\alpha}{4})$ which is  combined with Frostman's lemma.

For $s\in (1/2, d/2]$ let us define $\alpha_\ast(s)$  to be the infimum of $\alpha$ such that  the estimate \eqref{suptR} holds.   
Then, by a standard argument (see, for example,  \cite[Appendix B]{BBCR} or \cite[Section 17]{Mattila2}), it follows 
$ 
\mathcal D_d( s )   \le \, \alpha_\ast(s) .
$
Thus Theorem \ref{divergence} is a consequence of the following. 
\begin{thm}\label{thm:sup}
Let $d\ge3$ and  $\mu \in \mathfrak M^{d}(\alpha)$. 
Then we have \eqref{suptR} 
if  
\begin{align}\label{s22}
s> s(\alpha,d):=
\begin{cases}
\, \frac{d-\alpha}2, 
  & \  \hfill 0 < \alpha \le \frac{d-1}2, 
 \\
\, \frac{3d+1}8-\frac \alpha 4, 
 & \ 
\frac{d-1}2 < \alpha  \le \frac{d+1}2,  \\
\, \frac{d-\alpha}2+\frac{\alpha-1}{2(d-1)}, 
  & \  \hfill
\frac{d+1}2 < \alpha  \le d .  
\end{cases}
\end{align}
\end{thm}
Especially $s( \alpha,3)= \frac{5-\alpha}{4}$ for $\alpha\in  (1, 3] $, which is sharp.

Concerning the proof of Theorem \ref{thm:sup}, our key observation is that  the maximal estimates with respect to an $\alpha$-dimensional measure in $\mathbb R^d$ can be deduced from a family of estimates with respect to $\alpha$-dimensional measures in  $\mathbb R^{d+1}$. This can be shown by the Kolmogorov-Seliverstov-Plessner linearization argument. 
Thus  the Sobolev exponent $s$ in \eqref{suptR} is closely related to the exponent $\gamma=\gamma(\alpha)$ for which  the {\sl fractal Strichartz estimate}
\begin{equation}\label{fracS}
\| e^{it\sqrt{-\Delta}} f\|_{L^2(d\nu)}
\lesssim \langle\nu\rangle_\alpha^{\frac 12}  \|f\|_{H^\gamma(\mathbb R^d)}
\end{equation}  
holds  for any $\nu \in \mathfrak M^{d+1}(\alpha)$.   
For $d\ge2$, this type of estimates was studied in connection with geometric measure theory (e.g., \cite{Wolff3, Erdogan, O, OO, CHL, Harris2, Rogers}).  
In particular, the sharp exponent $\gamma$ was established  
by Wolff \cite{Wolff3},  Erdogan  \cite{Erdogan} when $d=2$  and by Cho and two  of the authors  \cite{CHL} when $d=3$. 
We discuss the estimate further  in Section 2 before the proof of Theorem \ref{divergence}.

\subsection{\texorpdfstring{$L^2$}{L2} estimate for \texorpdfstring{$e^{it \sqrt{-\Delta}}$}{} with a product measure}
The average decay estimate 
\eqref{scale-beta2}
has been of interest in relation to not only the Falconer distance set problem but also its refined variant: the pinned distance set problem.   
We refer the reader to  \cite{DZ, Liu, GIOW, Liu3, DIOWZ}  for the most recent progresses related to the distance set problems.

By the argument due to Mattila \cite{Mattila}, it can be shown that if $E$ is a Borel set, then  
\begin{equation}\label{dist}
\dim E >\alpha  \quad  \Rightarrow  \quad   |\{ |x-y| :x,y\in E \}| >0
\end{equation}  
provided that  \eqref{scale-beta2}
holds and  $\alpha +\beta >d$.

On the other hand,   Rogers \cite{Rogers} proposed an approach to Falconer's distance set problem which is based on,  instead of \eqref{scale-beta2},   a Strichartz type estimate for $\wa f$ with respect to fractal measures. Using the Riesz representation theorem, he showed that \eqref{dist} holds for $d-1 <2\alpha < d+1$ if the estimate
\begin{equation}\label{st}
\| \wa \pl  f \|_{L^2(\mathbb R^d\times I; d\mu dt)}  \lesssim \lambda^{\frac{d-1}{2} - \frac{\gamma}{2}} \|f\|_{L^2(\mathbb R^d)},   \quad  \mu\in  \mathfrak M^d(\alpha) 
\end{equation}
holds true for  $ \gamma > d - \alpha$.
Here, $I=(1,2)$ and $P_\lambda$ is the standard Littlewood-Paley projection operator which is given by  $\widehat{P_\lambda f}(\xi) =\psi(|\xi|/\lambda)\widehat f(\xi)$ where $\psi\in \mathrm C_c^\infty((2^{-1}, 2))$.

Interestingly the estimates  \eqref{scale-beta2} and \eqref{st} turn out to be essentially equivalent to each other.  
 
 \begin{thm}\label{equivlem2}
Let $\lambda \gg 1$ and $\mu \in \mathfrak M^d(\alpha)$ for $0< \alpha\le d$.  Then  \eqref{st}  holds for $\gamma<\beta$  if the estimate   \eqref{scale-beta2}  holds. 
Conversely, we have  the estimate   \eqref{scale-beta2}  with $\beta = \gamma $   if  \eqref{st}  holds. 
\end{thm}

Recently, Liu \cite{Liu} deduced the similar implication  as  \eqref{dist}  for the pinned distance set: there exists $x\in E$ such that $D_x(E)=\{|x-y|: y\in E\}$ has positive Lebesgue measure if the estimate  \eqref{scale-beta2}  holds and $\alpha + \beta >d$. The estimate 
\begin{equation}\label{sphe}
\| f\ast \sigma_t \|_{L^2(\mathbb B^d(0,1)\times I;d\mu dt)} \lesssim \|f\|_{H^{ - \frac{\beta}{2}+\epsilon} }    
\end{equation} 
played a crucial role in Liu's argument where $\sigma_t$ is the normalized surface measure on the sphere $\{x:|x|=t\}$ and $I=(1,2)$.  Using  an $L^2$ identity, he showed that  \eqref{scale-beta2}  implies  
\eqref{sphe} for any $\epsilon>0$.  
The estimate \eqref{sphe} is an easy consequence of \eqref{st} since 
$\widehat{d\sigma_t}=c_+ e^{it|\xi|}a_+(t\xi)+c_- e^{-it|\xi|}a_-(t\xi)$  with $a_\pm(t\xi)=O(|t\xi|^{-\frac{d-1}2})$ when $|\xi|\ge 1$. 
Using Theorem \ref{equivlem2} and the argument due to Mattila \cite{Mattila},  one can deduce  the results in Liu \cite{Liu} and Rogers \cite{Rogers}.

\section{Proof of Theorem \ref{thm:sup} and Theorem \ref{equivlem2}} \label{Sec:2}

To prove Theorem \ref{thm:sup},
we use the fractal Strichartz estimate \eqref{fracS}.
Concerning the Lebesgue measure the estimate \eqref{fracS} is well understood but for the general measure $\nu$ the sharp regularity is not known  when  $d\ge 4$. 
When $d=2$,  Wolff \cite{Wolff3} obtained the estimate \eqref{fracS}
with $\alpha <1$ and the sharp estimates for $\alpha\ge1$ was settled by Erdo\u{g}an \cite{Erdogan}. Those results were extended to higher dimensions ($d \ge 3$) by Cho and two of the authors    \cite{CHL}, where,  especially for $d=3$, the sharp estimate was proved by using the sharp bilinear restriction estimate for the cone.   
Recently, Harris \cite{Harris2}  improved the bound when $d\ge 4$ and $\frac{d+1}{2}<\alpha \le d$.  
These  results can be summarized as follows. 
\begin{thm}[\cite{CHL, Harris2}]\label{3d-frac}
Let $d\ge3$. For $0<\alpha < d+1$, let $\nu \in \mathfrak M^{d+1}(\alpha)$. 
Then \eqref{fracS}
holds  with $\gamma > \gamma(\alpha,d)$ where
\begin{align}\label{s2}
\gamma(\alpha,d)=
\begin{cases}
\frac{d-\alpha}2, 
\quad & \quad \hfill
0 < \alpha \le \frac{d-1}2,  \\
\frac{3d+1}8-\frac \alpha 4, 
\quad &  \quad \hfill
\frac{d-1}2 < \alpha  \le \frac{d+1}2,  \\
\frac{d-\alpha}2+\frac{\alpha-1}{2(d-1)}, 
\quad & \quad \hfill
\frac{d+1}2 < \alpha  \le d, \\
\frac{d+1-\alpha}2,
\quad & \quad \hfill
d < \alpha \le d+1.
\end{cases}
\end{align}
\end{thm}

Now we prove Theorem \ref{thm:sup} making use of Theorem \ref{3d-frac}.

\begin{proof}[Proof of Theorem \ref{thm:sup}]
In order to prove  Theorem \ref{thm:sup}, by Kolmogorov-Seliverstov-Plessner linearization argument,  it suffices to show that 
\begin{equation}\label{suptR1}
\|  e^{i\mathbf t(\cdot) \sqrt{-\Delta}}  f     \|_{L^2 (d\mu)} \leq C \langle \mu \rangle_\alpha^{\frac12} \|f\|_{H^s(\mathbb R^d)}
\end{equation}
for any measurable function $\mathbf t:\mathbb B^d(0,1)\to (0,1)$ with $C$ independent of $\mathbf t$. 

Let us define a linear functional $\ell$ by 
\[
\ell(F) = \int F(x, \mathbf t(x)) d\mu(x) 
\]
for any continuous function $F$ in $\mathrm C_0(\mathbb R^{d+1})$. 
By the Riesz representation theorem it follows that there is a unique Radon measure $\nu$ on $\mathbb R^{d+1}$ such that 
\[
\ell(F) = \int F(y,s) d\nu(y,s) = \int F(x, \mathbf t(x)) d\mu(x) . 
\]
Clearly $\nu$ belongs to $\mathfrak M^{d+1}(\alpha )$ and 
$\langle \nu \rangle_\alpha \le \langle \mu \rangle_\alpha$.\footnote{It is clear that $|(x-y,\mathbf t(x)-s)| \le r $ implies $|x-y| \le r$ and
\[\nu \big(\mathbb B^{d+1}((y,s),r) \big) =\int \chi_{\mathbb B^{d+1}((y,s),r) } \big(x,\mathbf t(x) \big)\,d\mu(x)
\le \int \chi_{\mathbb B^d(y,r)}(x)\,d\mu(x) \le \langle \mu\rangle_\alpha  r^\alpha. \]
} 
Since $\| e^{i\mathbf t(\cdot)\sqrt{-\Delta}}f \|_{L^2 ( d\mu)}  = \| \wa f\|_{L^2 (d\nu)}$, 
applying Theorem \ref{3d-frac}, we see that \eqref{suptR} holds if  \eqref{s22} is satisfied.
\end{proof}

We now proceed to prove Theorem \ref{equivlem2}. 
The relation of \eqref{scale-beta2}  and \eqref{st} is basically  due to the following lemma.

\begin{lem} \label{equivlem}
Let $q \ge 2$ and   $\lambda \gg 1$.
For $\mu \in \mathfrak M^d(\alpha)$, suppose that  
\begin{equation}\label{sp-2q}
\Big(\int |\widehat F|^q d\mu    \Big)^{\frac1q} \lesssim    \lambda^{s} 
\|F\|_2
\end{equation}
whenever $\supp F \subset \lambda \mathbb S^{d-1}+O(1)$. Then we have 
\begin{equation}\label{co-2q}
\Big( \int_1^2 \int | \widehat G |^q d\mu dt \Big)^{\frac1q} \lesssim 
\lambda^{s+\frac 12-\frac 1q}  
 \|G\|_2
\end{equation}
whenever $\supp G \subset \lambda \Gamma^d+O(1)$. Here $\Gamma^d:=\{ (\xi,|\xi|) \in \mathbb R^{d+1}: 2^{-1} \le |\xi| \le 2\}$.
In particular,  when $q=2$,  the estimates \eqref{sp-2q} and \eqref{co-2q} are equivalent. 
\end{lem}

\begin{proof}
Assume that \eqref{sp-2q} holds.
Let $\psi$ be a smooth function such that 
$\widehat \psi\gtrsim 1$ on $[1,2]$  and $\psi$ is supported in $[-1,1]$.   Then it follows 
\begin{equation}\label{Q1}
 \int_1^2 \int  \big|\widehat{G}(x,t) \big|^q  d\mu(x) dt \lesssim  \int \int_1^2   \big|\widehat{G} (x,t) \widehat \psi(t) \big|^q  dt d\mu(x). 
\end{equation}
We consider the inverse Fourier transform in $t$ and 
observe that $\big( \widehat G(x, \cdot) \widehat \psi(\cdot) \big)^\vee(\tau) =  \mathcal F_\xi G(x,\cdot) \ast \psi (\tau)$.
Here $\mathcal F_\xi$ denotes the Fourier transform in $\xi$, so that
 $\widehat G(x,t) = \int \mathcal F_\xi G(x,\tau) e^{i t\tau} d\tau$. 
For each $x$, it is easy to see that $\mathcal F_\xi G(x ,\cdot) \ast\psi$ is supported in $[\lambda /c,c\lambda ]$ for some $c >1$.
Thus, by Bernstein's inequality,  we obtain
\begin{equation}\label{Q2}
\int \int_1^2 \big|\widehat{G} (x,t) \widehat \psi(t) \big|^q  dt d\mu(x)
\lesssim \lambda^{(\frac 12-\frac 1q)q}
\int \Big( \int_1^2 \big| \widehat G(x,t) \widehat \psi(t) \big|^2\,dt \Big)^{\frac q2}\,d\mu(x).
\end{equation}
By Plancherel's theorem in $t$ and  Minkowski's inequality, we 
have 
\begin{equation}\label{Q3}
\begin{aligned}
 \int \! \Big(\! \int_1^2 \big| \widehat G(x,t) \widehat \psi(t) \big|^2\,dt \Big)^{\frac q2}\,d\mu (x)
 	= & 	\int\! \Big( \! \int  |\mathcal F_\xi G(x ,\cdot) \ast\psi(\tau) |^2 d\tau \Big)^{\frac q2} d\mu (x)  \\
  \lesssim  &
	\Big(\!\int\! \Big(\! \int  \big|\mathcal F_\xi G(x ,\cdot) \ast\psi(\tau) \big|^q d\mu (x)
	\Big)^{\frac 2q} d\tau \Big)^{\frac q2}.  
\end{aligned}
\end{equation}
Note that
\[
\mathcal F_\xi G(x ,\cdot) \ast\psi(\tau) =\mathcal F_\xi \B(\int {G} (\xi,\rho)  \psi(\tau-\rho) d\rho\B).
\]
Since $\psi$ is supported in $[-1,1]$ and $\supp G\subset   \lambda \Gamma^{d}+O(1)$, we see that for each fixed $\tau$,  $\int {G} (\xi,\rho)  \psi(\tau-\rho) d\rho$ 
 is supported in  the set $||\xi|-\tau|=O(1)$ and  $\tau\sim \lambda$. 
Thus, by \eqref{sp-2q} it follows that, for each $\tau$, 
\[   
\int  \big|\mathcal F_\xi G(x ,\cdot) \ast\psi(\tau) \big|^q  d\mu\lesssim     \lambda^{s q}\B\|\int {G} (\xi,\rho)  \psi(\tau-\rho) d\rho\B\|^q_{L^2_\xi}.
\]
Combining this with \eqref{Q1}, \eqref{Q2} and \eqref{Q3},  
  by Minkowski's and Young's convolution inequality we have 
\[ 
\Big( \int_1^2 \int  \big|\widehat{G} \big|^q  d\mu dt \Big)^{\frac 1q}
\lesssim    \lambda^{s+\frac 12-\frac 1q} 
\Big(\int \B\|\int {G} (\xi,\rho)  \psi(\tau-\rho) d\rho\B\|^2_{L^2_\xi} d\tau \Big)^{\frac 12} \lesssim    \lambda^{s+\frac 12-\frac 1q} \|G\|_2 .\]
Therefore we get \eqref{co-2q}.

Now let $q =2$.  
Let $F$ be a function with $\supp F\subset   \lambda \mathbb S^{d-1}+O(1)$ and consider 
\[ G(\xi,\tau)=F(\xi) \chi_{[\lambda, \lambda+1]}(\tau) . \] 
Then, since  $ \int_1^2|  \widehat{\chi_{[\lambda, \lambda+1]}}(t) |^2dt\gtrsim 1$ independent of $\lambda$,  by Fubini's theorem it is clear 
\[  \int  \big|\widehat{F} \big|^2 d\mu  \lesssim    \int_1^2 \int  \big|\widehat{G} \big|^2  d\mu dt . \]   
By \eqref{co-2q} with $q=2$ and Plancherel's theorem,  it follows that 
\[  \int  \big|\widehat{F} \big|^2 d\mu  \lesssim    \lambda^{2s}  \|G\|^2_2 =  \lambda^{2s} \|F\|_2^2.   \]
Hence \eqref{sp-2q} and \eqref{co-2q} are equivalent when $q=2$.  
 \end{proof}

Using this lemma, we prove Theorem \ref{equivlem2}.
\begin{proof}[Proof of Theorem \ref{equivlem2}]
Since $\mu$ has compact support,  by duality and the uncertainty principle, one can easily see that  \eqref{scale-beta2} is equivalent to 
\begin{equation}\label{sphere1} 
   \big(\int | \widehat{F} | d\mu \big)^2  \lesssim   \lambda^{d-1-\beta}  \|\mu\|\langle \mu\rangle_\alpha  \|F\|^2_2,
\end{equation}
when $\supp F \subset \lambda \mathbb S^{d-1}+O(1)$ and $\mu \in \mathfrak M^d(\alpha)$.
By the argument in  \cite{BBCR}
the estimate can be strengthened to  
\begin{equation}\label{weak22}
\|  \widehat F  \|_{L^{2,\infty}(d\nu)}\lesssim   \lambda^{(d-1-\beta)/2}  \langle \nu\rangle_\alpha^{1/2}   \|F\|_2
\end{equation} for $\nu \in \mathfrak M^d(\alpha)$. To see this, we consider a measure $d \nu_E = \nu(E)^{-1} \chi_E d\nu$ for $\nu \in \mathfrak M^d(\alpha)$ and a Borel set $E$. Then $\nu_E \in \mathfrak M^d(\alpha)$, $\|\nu_E\|\le1$, and $\langle \nu_E \rangle_\alpha \le \nu(E)^{-1} \langle \nu \rangle_\alpha$. Applying \eqref{sphere1} with $\nu_E$, we see that $  {\nu(E)}^{-1/2}   \int_E |\widehat F|  d\nu   = \nu(E)^{1/2}   \int_E | \widehat F | d\nu_E    \lesssim \lambda^\frac{d-1-\beta}{2} \langle\nu\rangle_\alpha^{\frac12} \|F\|_2$. So, we get \eqref{weak22} by taking $E= \{ x: |\widehat F| > \omega \}$.     

Interpolating \eqref{weak22} with the trivial estimate $\|  \widehat F   \|_{L^{\infty}(d\nu)}\lesssim \lambda^{\frac{d-1}2}  \|F\|_2$ we obtain $L^2\to L^q(d\nu)$ estimate. Since $\supp \nu\subset \mathbb B^d(0,1)$, taking $q$ arbitrarily close to $2$, 
we get 
\begin{align}\label{sphere}
 \int  \big|\widehat{F} \big|^2 d\nu  \lesssim  \lambda^{d-1-\gamma}   \langle \nu\rangle_\alpha    \|F\|^2_2
\end{align}
for  $\gamma<\beta$.  
By Lemma \ref{equivlem}, we can see that \eqref{sphere} is equivalent to 
\begin{align}\label{cone}
\int_1^2 \int  \big|\widehat{G} \big|^2  d\nu dt \lesssim  \lambda^{d-1-\gamma}  \langle \nu\rangle_\alpha \|G\|^2_2,  
\end{align}
where $\supp G \subset   \lambda \Gamma^d+O(1)$. 
Let $\varphi\in \mathcal S(\mathbb R^d)$ such that $\varphi\ge 1$ on $\mathbb B^d(0,1)$ and $\supp \widehat \varphi\subset \mathbb B^d(0,1)$. 
Then the space time Fourier transform of $\varphi \wa P_\lambda f$ is supported in $\lambda\Gamma^d+O(1)$. By  \eqref{cone} and Plancherel's theorem 
we get
\begin{align}\label{ext-cone}
\int_1^2 \int  \B|\varphi (x)\wa P_\lambda f(x) \B|^2   d\nu dt\lesssim  \lambda^{d-1-\gamma}  \langle \nu\rangle_\alpha \|f\|_2^2,
\end{align}
which gives \eqref{st} with $\gamma <\beta $ as desired.

To show the converse implication we note that \eqref{st}  implies \eqref{cone}. 
It can be seen easily by Plancherel's theorem and the standard slicing argument decomposing the conic neighborhood into a family of cones.
Thus  by  Lemma \ref{equivlem} we have \eqref{sphere} which clearly implies 
$\big|\int  \widehat{F} d\mu \big|^2 \lesssim   \lambda^{d-1-\gamma} \|\mu\|\langle \mu\rangle_\alpha \|F\|_2^2$. Therefore by duality and the uncertainty principle we get  \eqref{scale-beta2}  with $\beta=\gamma$. 
\end{proof}

\appendix

\section{}

In this section, we discuss necessary conditions for the  estimate
\begin{equation}\label{suptR-q}
\| \sup_{0<t <1}  |\wa  f   |  \|_{L^q (d\mu)} \lesssim  \langle \mu \rangle_\alpha^{\frac1q} \|f\|_{ H^s(\mathbb R^d)} .
\end{equation}

\begin{lem}\label{nece}
Let $d\ge 3$. The estimate \eqref{suptR-q} 
holds for all $\mu \in \mathfrak M^d(\alpha)$  only if 
\begin{equation}\label{low-sup2}
s\ge 
\begin{cases}
\  \max(\frac d2-\frac \alpha q, \frac{d+1}4), & ~ 0<\alpha \le 1, 
\\ 
\ \max(\frac d2-\frac \alpha q, \frac{d+1}4-\frac{\alpha-1}{2q},
\frac{d+2-\alpha}4, \frac{d-\alpha}2), &~ 1<\alpha \le d.
\end{cases}
\end{equation}
\end{lem}

\begin{proof}[Proof of Lemma \ref{nece}]
We first show $s\ge \frac d2-\frac \alpha q$.
Let $\widehat f(\xi)=\chi_{\mathbb B^d(0,\lambda)}(\xi)$. If we take $t=1/\lambda$, then
we have $|e^{it\sqrt{-\Delta}}f(x)| \gtrsim \lambda^d$ for $|x| \le \lambda^{-1}$,
which implies
$
\sup_{0<t<1}|e^{it\sqrt{-\Delta}}f(x)| \ge \lambda^d
$
for $|x| \le \lambda^{-1}$. 
We take $d\mu=\chi_{\mathbb B^d(0,1/2)}(x)|x|^{\alpha-d}dx$. Then it is clear that   $\langle \mu \rangle_\alpha \lesssim 1$. 
A simple computation gives 
$\lambda^d \lambda^{-\frac \alpha q} \lesssim \lambda^{s+\frac d2}$.
This shows $s \ge \frac d2-\frac \alpha q$.

Next, we show that 
\begin{equation}\label{sss}
s \ge \begin{cases}
\, \frac{d+1}{4}, & ~ 0<\alpha \le 1,\\
\, \frac{d+1}4-\frac{\alpha-1}{2q}, & ~ 1< \alpha \le d.
\end{cases}
\end{equation}
Let $P = \{ \xi = (\xi_1,\xi')\in \mathbb R\times \mathbb R^{d-1}: \lambda \le  \xi_1  \le 2\lambda, |\xi'|\le \lambda^{1/2} \}$. 
Then we consider  $f$ which is given by $\widehat f = \chi_P$. 
Then  the set $\{ (\xi,|\xi|): \xi\in P\}$ is contained in a $ (d+1) $-dimensional rectangle $Q$ of dimensions $C\lambda\times \underbrace{C\lambda^{1/2}\times \cdots \times C\lambda^{1/2}}_{d-1 \, times}\times C$ for a constant $C>0$. 
 We see that $ | \wa f | \gtrsim |P| $ on the set $Q^* = \{ (x_1, x', t) : |x_1 - t| \le c \lambda^{-1},  |x'| \le c \lambda^{-1/2},  t\le c \}$ for a constant $c \le 1$.  
It follows that $ \| \sup_{0<t<1} |\wa f(x)| \|_{L^q(\mu)} \gtrsim |P| \mu(Proj_x(Q^*))^{1/q}$.
Here, $Proj_x :\mathbb R^{d+1} \rightarrow \mathbb R^d$ is the orthogonal projection which maps $(x,t)$ to $(x,0)$.

For a given $\alpha$, let $\ell \in [1,d]$ be an integer such that $\ell-1 <\alpha\le \ell$. 
For a fixed $\ell$, we define a measure 
\[
d\mu  = \prod_{i=1}^{d-\ell} d\delta(x_{d-i+1}) |x_\ell|^{\alpha-\ell}dx_\ell dx_{\ell-1} \cdots d x_1,
\]
for the delta measure $\delta$ and $x = (x_1,\dots, x_d)$. 
Clearly, $\mu \in \mathfrak M^d(\alpha)$.

Since  $Proj_x(Q^*)$ is comparable to a rectangle given by $\{ (x_1,x')\in \mathbb R^d :   x_1 \lesssim 1, |x'| \le c \lambda^{-1/2} \}$ for a sufficiently large $\lambda$, it is easy to show that 
$\mu(Proj_x(Q^*)) \gtrsim 1$ for $0<\alpha\le 1$ or $\mu(Proj_x(Q^*)) \gtrsim \lambda^{-(\alpha-1)/2}$ for $1<\alpha \le d$. 
Thus we get $|P| \mu(Proj_x(Q^*) )^{\frac1q} \lesssim \lambda^s |P|^{\frac12}$ from \eqref{suptR-q}. This implies \eqref{sss}.

The proof of $s \ge \frac{d+2-\alpha}4$ can be shown similarly as  in \cite{CHL}, so we shall be brief.
Let $P$, $Q$, and $Q^*$ as above.
Also let $\phi_P$ be a Schwartz function supported in $P$. 
For an integer $N \sim \lambda^{\frac{\alpha-1}2}$, we take
\[
\widehat f(\xi) = N^{-\frac 12}\sum_{k=1}^N  \phi_P(\xi) e^{-iv_k \cdot \xi'} ,
\]
where $v_k$ are the lattice points on $\mathbb B^{d-1}(0,1)$ separated by $\sim \lambda^{-\frac{\alpha-1}{2(d-1)}}$.
Then we have  $ |e^{it\sqrt{-\Delta}}f(x)| \gtrsim N^{-\frac12} |P| \chi_S(x,t) $ for $S = \cup_{k=1}^N  (  Q^*  +(0,v_k,0)  )$.
Let $d\mu=\lambda^{\frac{d-\alpha}2}\chi_{Proj_x(S)} dx$, then $\langle \mu \rangle_\alpha \lesssim 1$ and $\mu(Proj_x(S))\sim1$.
Also,  $\|f\|_2^2 \lesssim \sum_k \| \widehat {f_k}\|_2^2 \lesssim |P| \sim \lambda^{\frac{d+1}2}$.
Therefore, 
we obtain
$N^{-\frac 12}\lambda^{\frac{d+1}2} \lesssim \|\sup_{0<t<1} |e^{it\sqrt{-\Delta}}f| \|_{L^q(\mu)} \lesssim \|f\|_{H^s} \lesssim \lambda^s \lambda^{\frac{d+1}4}$.
Since $N\sim \lambda^{\frac{\alpha-1}2}$, we get
$s \ge \frac{d+2-\alpha}4$.

In a similar way,  we can prove $s \ge \frac{d-\alpha}2$. (See also \cite{Rogers} for an alternative proof.)
Let $M$ be an integer satisfying $M \sim \lambda^\alpha$.   
Let us set
\[
\widehat f(\xi)= 
M^{-\frac 12}\sum_{k=1}^M  \phi_{\mathbb B^d(0,\lambda)}(\xi)e^{-i\omega_k \cdot \xi} ,
\]
where $\omega_k$ are lattice points in $\mathbb B^d(0,1)$ separated by $\lambda^{-\frac \alpha d}$.
If we choose $t=1/\lambda$, then
$\sup_{0<t<1} |e^{it\sqrt{-\Delta}} f (x,t)| \ge M^{-\frac 12}\lambda^d \chi_T(x) $, where $T:=\cup_{k=1}^M \mathbb B^d(\omega_k,\lambda^{-1})$.
Let $d\mu:=\lambda^{d-\alpha}\chi_T(x)\,dx$, then
$\langle \mu \rangle_\alpha \lesssim 1$ and $\mu(T) \sim 1$.
We can conclude that
$M^{-\frac 12}\lambda^{d} \lesssim \|\sup_{0<t<1} |e^{it\sqrt{-\Delta}} f| \|_{L^q(\mu)} \lesssim \|f\|_{H^s} \lesssim \lambda^s \lambda^{\frac d2}$, which gives
$s \ge \frac{d-\alpha}2$.  
\end{proof}

The following shows failure of weak type estimate when $s= 1/2$. 

\begin{lem}\label{weak-type}   Let $q < \infty$ and let $\mu\in \mathfrak M^d(\alpha)$. Suppose $\|\mu\|\neq 0$, then the estimate  
$\| \sup_{0<t <1}  |\wa  f   |  \|_{L^{q,\infty} (d\mu)} \lesssim  \langle \mu \rangle_\alpha^{\frac1q} \|f\|_{ H^\frac12(\mathbb R^d)}$ fails.
\end{lem} 

\begin{proof}
Since $\|\mu\|\neq 0$, there is an $r>0$ such that $\mu(\{x: r\le |x|\le 1\} )\neq 0$. 
Let $\beta\in C_c((1/2,2))$ such that $\int \beta(s)  s^{\frac{d-1}{2}}  ds=1/ (2\pi)^{\frac d2} $. Setting $\beta_N(\xi)=N^{-\frac{d+1}2} \beta(|\xi|/N)$ for $N\ge  1/r$, we  consider 
\[  
\int e^{i(x\cdot \xi + t|\xi|)} N^{-\frac{d+1}2}  \beta(|\xi|/N)  d\xi 
=  N^{\frac{d-1}{2}}  \int  \int_{\mathbb S^{d-1}} e^{  iNs x\cdot \theta} d\theta e^{ iNts} \beta(s) s^{d-1} ds.
\]
Since $\int_{\mathbb S^{d-1}} e^{iNs x\cdot \theta} d\theta=  (2\pi)^{\frac d2}  (Ns|x|)^\frac{2-d}2 J_{\frac{d-2}2} (Ns|x|)$, using asymptotic expansion of the Bessel function 
($J_{\frac{d-2}2} (s)=C_\pm s^{-\frac12}e^{\pm is} +O(s^{-\frac32}), s\ge 1),$
and taking $t=|x|$, 
we see 
 \Be 
 \label{asym}
 \int e^{i(x\cdot \xi + |x||\xi|)} \beta_N(\xi)  d\xi=  C_{-}|x|^{-\frac{d-1}2} +O(N^{-1}), \quad r\le |x|\le 1
 \Ee 
 provided that $N\ge 1/r$. 
Let us consider $f_L$ which is given by $ \widehat{f_L}=  \sum_{1/r\le 2^k\le L} \beta_{2^k}$.
Then by \eqref{asym} 
we see that 
\[  e^{i|x| \sqrt{-\Delta}}f_L(x)=  \sum_{1/r\le 2^k\le L} C_{-}|x|^{-\frac{d-1}2}   +O(1).  \] 
Thus, $\sup_{0<t <1}  |\wa  f_L (x) |\sim \log L$ if $r\le |x|\le 1$ and $L\ge 10^2/r^2$. If we assume the estimate 
$\| \sup_{0<t <1}  |\wa  f   |  \|_{L^{q,\infty} (d\mu)} \lesssim  \langle \mu \rangle_\alpha^{\frac1q} \|f\|_{ H^s(\mathbb R^d)}$, then 
  it follows that  
$(\log L) \mu(\{x:   r  \le |x|\le 1\} )^\frac{1}{q} \lesssim  \langle \mu \rangle_\alpha^{\frac1q}(\log L)^\frac12$ 
because $\| f_L \|_{H^\frac12}\sim (\log L)^\frac12$.  This can not be true unless $\mu(\{x: r\le |x|\le 1\} )=0$. 
This is a contradiction. 
\end{proof}

\subsection*{Acknowledgement} 
This work was supported by the National Research Foundation of Korea (NRF) grant number NRF-2017R1C1B2002959 (Seheon Ham), NRF-2019R1A6A3A01092525 (Hyerim Ko), and NRF-2018R1A2B2006298 (Sanghyuk Lee).

\bibliographystyle{plain}

\end{document}